\documentclass{article}
\usepackage{amstex,amssymb}
\begin{document}
\newenvironment{rlist}{\begin{list}{$(\roman{enumi})$}
{\usecounter{enumi}
\setlength{\rightmargin}{10pt}
\setlength{\leftmargin}{40pt}
\setlength{\labelwidth}{30pt}
\setlength{\parsep}{1pt}
\setlength{\itemsep}{1pt}}}{\end{list}}
\newtheorem{thm}{Theorem}[section]
\def\Aut{\mathop{\rm Aut}}
\def\Hol{{\textstyle\mathop{\rm Hol}}}
\def\eq#1{(\ref{#1})}
\def\R{\mathbin{\mathbb R}}
\def\Z{\mathbin{\mathbb Z}}
\def\C{\mathbin{\mathbb C}}
\def\CP{\mathbin{\mathbb{CP}}}
\def\ha{{\textstyle{1\over2}}}
\def\md#1{\vert#1\vert}

\title{On the topology of desingularizations of Calabi-Yau orbifolds}
\author{Dominic Joyce \\
Lincoln College, Oxford, OX1 3DR, England}
\date{June 1998}
\maketitle

\section{Introduction}
A {\it Calabi-Yau 3-fold} is a compact complex 3-manifold $(X,J)$
equipped with a K\"ahler metric $g$ and a holomorphic volume form 
$\Omega$, which is a nonzero $(3,0)$-form on $X$ with $\nabla\Omega=0$, 
where $\nabla$ is the Levi-Civita connection of $g$. Calabi-Yau 3-folds 
are Ricci-flat and have $\Hol(g)\subseteq SU(3)$, and one can construct
many examples of them using Yau's solution of the Calabi conjecture.

Suppose $X$ is a Calabi-Yau 3-fold and $G$ a finite group that acts 
on $X$ preserving $J,g$ and $\Omega$. Then $X/G$ is a Calabi-Yau 
orbifold, where an orbifold is a singular manifold with at most 
quotient singularities. Often it is possible to find a compact 
manifold $Y$ that desingularizes $X/G$, and carries a family of 
Calabi-Yau structures that converge to the singular Calabi-Yau 
structure on $X/G$ in a well-defined sense, so that the orbifold 
metric on $X/G$ may be regarded as the degenerate case in a smooth 
family of Calabi-Yau metrics on~$Y$. 

In this paper we will study Calabi-Yau 3-folds $Y$ that arise by 
deforming the complex structure of orbifolds $X/G$ which have
singularities in complex codimension two. We shall focus on the 
topology of $Y$, and will show that a single orbifold $X/G$ can
admit several different Calabi-Yau desingularizations $Y_1,Y_2,
\dots,Y_k$, with different Hodge numbers and Euler characteristics. 
We give an explanation of this phenomenon in terms of the `Weyl 
group' of the codimension two singularities, and make some tentative
suggestions on how this idea may help to explain our examples in the 
physical language of String Theory. We also briefly discuss
resolutions of 7-and 8-orbifolds with the exceptional holonomy
groups $G_2$ and~$Spin(7)$.

The two main strategies for desingularizing $X/G$ to get $Y$ are 
called {\it resolution} and {\it deformation}. A {\it resolution} 
$(Y,\pi)$ of $X/G$ is a nonsingular complex 3-fold $Y$ with a proper 
holomorphic map $\pi:Y\rightarrow X/G$, that induces a biholomorphism 
between dense open subsets of $Y$ and $X/G$. We call $Y$ a {\it 
crepant resolution} if $K_Y\cong\pi^*(K_{X/G})$. Since Calabi-Yau 
manifolds have trivial canonical bundles, to get a Calabi-Yau 
structure on $Y$ we must choose a crepant resolution. 

A {\it family of deformations} of $X/G$ consists of a (singular)
complex 4-manifold $\cal Y$ with a proper holomorphic map 
$f:{\cal Y}\rightarrow\Delta$, where $\Delta$ is the unit disc 
in $\C$, such that $Y_0=f^{-1}(0)$ is isomorphic to $X/G$. The other 
fibres $Y_t=f^{-1}(t)$ for $t\ne 0$ are called {\it deformations} of 
$X/G$. If the $Y_t$ are smooth manifolds for $t\ne 0$, they are called 
{\it smoothings} of $X/G$. 

Thus, resolutions and smoothings are two different ways to 
desingularize a complex orbifold. We can also combine the two by 
smoothing a partial resolution of $Z$, or by resolving a deformation 
of $Z$. We shall use the word {\it desingularization} to mean any of 
these processes.

A number of important ideas about the topology of desingularizations 
$Y$ of a Calabi-Yau orbifold $X/G$ were first proposed by physicists 
working in String Theory, motivated by physical considerations. In 
1985, Dixon et al.~\cite[p.~684]{DHVW} conjectured that the Euler 
characteristic of $Y$ should be given by
\begin{equation}
\chi(Y)=\chi(X,G)={1\over\md{G}}\sum_{g,h\in G:\atop gh=hg}
\chi(X^g\cap X^h),
\label{peuler}
\end{equation}
where $X^g$ is the submanifold of $X$ fixed by $g\in G$, and $\chi$ 
is the Euler characteristic. Vafa \cite{Vafa1} and Zaslow 
\cite[p.~312]{Zasl} (see also Batyrev and Dais \cite{BaDa}) 
found a related formula for the Hodge numbers $h^{p,q}(Y)$ of $Y$, 
which we will not give.

As far as the author understands, the physicists who made these
conjectures believed that their formulae should apply to {\it at 
least one} desingularization $Y$ of $X/G$, but not necessarily to 
every desingularization; and an example of a desingularization
$Y$ for which \eq{peuler} does not hold appears in a paper by
some of the same authors~\cite[\S 2]{VaWi}.

Much more is known about the case of {\it crepant resolutions\/} 
$Y$ of $X/G$, particularly in complex dimension three. Roan 
\cite{Roan} proves that every 3-dimensional Calabi-Yau orbifold
admits a crepant resolution, for which \eq{peuler} holds. An 
orbifold $X/G$ can admit several topologically distinct crepant
resolutions $Y_1,\dots,Y_k$. But these resolutions are all related
by `flops', and have the same Euler characteristic and Hodge numbers.

Reid \cite{Reid} and Ito and Reid \cite{ItRe} develop a theory of
crepant resolutions of orbifolds which they call the `McKay 
correspondence', and in dimension 3 they show that every crepant
resolution $Y$ of $X/G$ must satisfy \eq{peuler}, and also
the Hodge number formulae of Vafa and Zaslow. In dimensions
four and above, it is known that some orbifolds admit no crepant
resolution (since they have `terminal singularities'). However,
there is strong evidence that crepant resolutions in all dimensions
must satisfy \eq{peuler} and the Hodge number formulae.

These results still leave open the question of the topology of
Calabi-Yau desingularizations $Y$ of $X/G$ that are not crepant
resolutions, but involve some deformation of the complex structure.
What can we say about this situation? By Schlessinger's Rigidity 
Theorem \cite{Schl}, quotient singularities of codimension 3 or 
more have no nontrivial deformations. But Calabi-Yau orbifolds
cannot have singularities in codimension one. So to resolve
$X/G$ by deformation, the singularities must be of {\it 
codimension two}.

The simplest sort of codimension two singularities in Calabi-Yau
3-folds are modelled on $\C\times(\C^2/H)$, where $H$ is a 
finite subgroup of $SU(2)$. The natural way to desingularize this 
is to use $\C\times X$, where $X$ is a Calabi-Yau desingularization
of $\C^2/H$. In fact, Kronheimer \cite{Kron1,Kron2} shows that we 
can give $X$ a metric with holonomy $SU(2)$ that is asymptotic to 
the flat metric on $\C^2/H$ at infinity, making $X$ into an {\it 
ALE space}. However, ALE spaces are all diffeomorphic to the
(unique) crepant resolution of $\C^2/H$, so in this case too, any
desingularization has the topology of a crepant resolution.

Therefore, to find desingularizations $Y$ which do not have the
topology of crepant resolutions, we must consider singularities
modelled on $\C^3/G$, where $G$ is a finite subgroup of $SU(3)$ 
with a nontrivial subgroup $H$ contained in some $SU(2)\subset 
SU(3)$, to produce the codimension two singularities, but where
$G$ is not itself contained in $SU(2)$. For instance, consider
a finite subgroup $G$ of elements of $SU(3)$ of the form
\begin{equation}
\begin{pmatrix}e^{i\theta}&0&0\\0&a&b\\0&c&d\end{pmatrix},
\quad\text{where}\quad\begin{pmatrix}a&b\\c&d\end{pmatrix}\in
U(2)\quad\text{and}\quad ad-bc=e^{-i\theta}.
\end{equation}
Let $H$ be the subgroup of elements $h\in G$ for which $e^{i\theta}=1$.
Then $H$ is a finite subgroup of $SU(2)$. It is also easy to
show that $H$ is a normal subgroup of $G$ and that $G/H\cong\Z_k$
for some positive integer $k$, where $e^{i\theta}$ is a $k^{\rm th}$
root of unity for each~$g\in G$. 

To desingularize $\C^3/G$, we may proceed in two stages. The first
stage is to choose a desingularization $X$ of $\C^2/H$, so that 
$\C\times X$ is a resolution of $\C^3/H$. Then we hope to find an 
action of $\Z_k=G/H$ upon $\C\times X$, which is asymptotic to the
prescribed action of $\Z_k$ on $\C^3/H$. If we can find such an
action, then the second stage is to desingularize $(\C\times X)/\Z_k$, 
either by a crepant resolution or a smoothing, to get a 
desingularization $Y$ of~$\C^3/G$.

Our key observation is the following: although the diffeomorphism 
type of $X$ is uniquely determined by $H$, the action of $\Z_k$ on 
$X$ and on its cohomology is not always uniquely determined by $G$. 
Instead, there can be a finite number of topologically distinct 
ways for $\Z_k$ to act on $X$, depending on the choice of complex 
structure of $X$, and on the level of cohomology these actions 
differ by an element of the {\it Weyl group} of the singularity 
$\C^2/H$. For one of these $\Z_k$-actions the desingularization 
$Y$ of $(\C\times X)/\Z_k$ has the topology of a crepant resolution, 
but for other choices of the $\Z_k$-action $Y$ does not have this 
topology, and its Euler characteristic is not given by~\eq{peuler}.

The rest of the paper is set out as follows. Section 2 summarizes 
the theory of singularities $\C^2/H$ for $H$ a finite subgroup of 
$SU(2)$, their desingularizations, and the idea of the Weyl group. 
Section 3 explains some theory about the orbifolds we are interested 
in, and the topology of their desingularizations, and \S 4 gives an 
example of this. Section 5 discusses a second way in which an orbifold
can have several Calabi-Yau desingularizations with different topology, 
and \S 6 considers another example, the orbifold $T^6/\Z_2^2$, which 
combines both phenomena, and has other interesting features.

In \S 7 we extend the discussion to the exceptional holonomy groups 
$G_2$ and $Spin(7)$, and give an example of an isolated singularity 
$\R^8/G$ which fits into our framework. Finally, section 8 speculates 
about the interpretation of these examples in String Theory. The 
author would like to thank Cumrun Vafa and David Morrison for helpful 
suggestions, and for correcting some of his misunderstandings 
about String Theory.

\section{Kleinian singularities and ALE spaces}

The quotient singularities $\C^2/H$, for $H$ a finite subgroup 
of $SU(2)$, were first classified by Klein in 1884 and are
called {\it Kleinian singularities;} they are also called 
{\it Du Val surface singularities}, or {\it rational double 
points}. The theory of these singularities and their resolutions
is very rich, and has many connections to other areas of
mathematics. Most of the following facts are taken from McKay 
\cite{McKa}, Slodowy \cite{Slod}, and Kronheimer \cite{Kron1,Kron2}.
A good reference on Lie groups, Dynkin diagrams and Weyl groups
is Bourbaki \cite{Bour}, in particular the tables on pages~250-270.

There is a 1-1 correspondence between finite subgroups $H\subset 
SU(2)$ and the {\it Dynkin diagrams} of type $A_r$ ($r\ge 0$), 
$D_r$ ($r\ge 4$), $E_6$, $E_7$ and $E_8$. Let $\Gamma$ be the
Dynkin diagram associated to $H$. These Dynkin diagrams appear 
in the classification of Lie groups, and each one corresponds 
to a unique compact, simple Lie group; they are the set of such 
diagrams containing no double or triple edges.

Each singularity $\C^2/H$ admits a unique crepant resolution
$(X,\pi)$. The preimage $\pi^{-1}(0)$ of the singular point 
is a union of a finite number of {\it rational curves} in $X$.
These curves correspond naturally to the vertices of $\Gamma$.
They all have self-intersection $-2$, and two curves intersect 
transversely at one point if and only if the corresponding 
vertices are joined by an edge in the diagram; otherwise the 
curves do not intersect. 

These curves give a basis for the homology group $H_2(X,\Z)$, 
which may be identified with the {\it root lattice} of the diagram, 
and the intersection form with respect to this basis is the negative
of the Cartan matrix of $\Gamma$. Define $\Delta$ to be 
$\{\delta\in H_2(X,\Z):\delta\cdot\delta=-2\}$. Then $\Delta$ 
is identified with the {\it set of roots} of the diagram. There are 
also 1-1 correspondences between the curves and the nonidentity 
conjugacy classes in $H$, and also the nontrivial representations 
of $H$; it makes sense to regard the nonidentity conjugacy classes 
as a basis for $H_2(X,\Z)$, and the nontrivial representations as 
a basis for~$H^2(X,\Z)$.

By the theory of Lie groups, the Dynkin diagram $\Gamma$ of 
$\C^2/H$ has a {\it Weyl group} $W$, and a representation of 
$W$ on the root lattice $H_2(X,\Z)$ of $\Gamma$. This action 
of $W$ preserves the subset $\Delta$ and the intersection form 
on $H_2(X,\Z)$, and by duality $W$ also acts on $H^2(X,\Z)$. 
Let $\Aut(\Gamma)$ be the group of automorphisms of the graph 
$\Gamma$, which is given by
\begin{equation}
\Aut(\Gamma)=\begin{cases}
\{1\} & \text{if $\Gamma=A_1,E_7$ or $E_8$,} \\
\Z_2 & \text{if $\Gamma=A_k$ ($k\ge 2$), 
$D_k$ ($k\ge 5$) or $E_6$,} \\
S_3 & \text{if $\Gamma=D_4$.} \end{cases}
\end{equation}

Now the vertices of $\Gamma$ correspond to the basis elements of
$H_2(X,\Z)$, so that $\Aut(\Gamma)$ acts naturally on $H_2(X,\Z)$,
preserving the intersection form. But the Weyl group $W$ also acts 
on $H_2(X,\Z)$. It turns out that there is a natural semidirect
product $\Aut(\Gamma)\ltimes W$, and the actions of $\Aut(\Gamma)$
and $W$ on $H_2(X,\Z)$ combine to give a representation of 
$\Aut(\Gamma)\ltimes W$ on $H_2(X,\Z)$. Define $\rho$ to be the 
dual representation of $\Aut(\Gamma)\ltimes W$ on both $H^2(X,\R)$
and $H^2(X,\C)$. The action of each element of $\Aut(\Gamma)\ltimes W$ 
is induced by a diffeomorphism of $X$, so we can interpret 
$\Aut(\Gamma)\ltimes W$ as a group of {\it isotopy classes of 
diffeomorphisms} of $X$. However, in general $\Aut(\Gamma)\ltimes W$ is 
not a group of diffeomorphisms of $X$, nor an isometry group of 
any of the metrics or complex structures on~$X$.

The singularities $\C^2/H$ can be desingularized by {\it deformation}
as well as by crepant resolution. Klein found that each singularity 
$\C^2/H$ is isomorphic as an affine complex variety to the zeros of 
a polynomial on $\C^3$. For example, $\C^2/\Z_k$ may be identified 
with the set of points $(x,y,z)\in\C^3$ for which $xy-z^k=0$. The 
deformations of $\C^2/H$ are constructed by adding terms of lower 
order in $x,y$ and $z$ to this polynomial. All of the smooth 
deformations of $\C^2/H$ are diffeomorphic to the unique crepant 
resolution $X$ of~$\C^2/H$.

Each of these desingularizations of $\C^2/H$ carries a special 
family of metrics with holonomy $SU(2)$. The metrics are asymptotic
up to $O(r^{-4})$ to the Euclidean metric on $\C^2/H$, and so are 
called {\it Asymptotically Locally Euclidean}; the complex manifold 
$X$ with its K\"ahler metric is called an {\it ALE space}. A complete 
construction and classification of ALE spaces was carried out by 
Kronheimer \cite{Kron1,Kron2}, and we describe it next.

Let $Y$ be an ALE space asymptotic to $\C^2/H$. Then $Y$ is 
diffeomorphic to $X$, and carries a geometric structure which
is encoded in the {\it K\"ahler form} $\omega$ and the {\it holomorphic 
volume form} $\Omega$ of $X$. Both $\omega$ and $\Omega$ are closed
forms, so they define de Rham cohomology classes $\alpha=[\omega]\in
H^2(X,\R)$ and $\beta=[\Omega]\in H^2(X,\C)$. Thus, to each ALE 
space $Y$ we may associate the pair~$(\alpha,\beta)\in H^2(X,\R)\times 
H^2(X,\C)$.

For each pair $(\alpha,\beta)\in H^2(X,\R)\times H^2(X,\C)$, Kronheimer 
\cite{Kron1} defined an explicit, possibly singular ALE space 
$X_{\alpha,\beta}$ asymptotic to $\C^2/H$, using the hyperk\"ahler 
quotient construction. Let $U$ be the subset
\begin{equation}
\!\!\!\!\!\!
U\!=\!\bigl\{(\alpha,\beta)\!\in\!H^2(X,\R)\!\times\!H^2(X,\C)\!:\!
\text{$\alpha(\delta)\!\ne\!0$ or $\beta(\delta)\!\ne\!0$
for all $\delta\in\Delta$}\bigr\}.
\end{equation}
Then $U$ is a dense open subset of $H^2(X,\R)\times H^2(X,\C)$.
Kronheimer showed that if $(\alpha,\beta)\notin U$ then 
$X_{\alpha,\beta}$ is an orbifold, and if $(\alpha,\beta)\in U$
then $X_{\alpha,\beta}$ is nonsingular and diffeomorphic to $X$,
and the K\"ahler form $\omega$ and holomorphic volume form $\Omega$ of
$X_{\alpha,\beta}$ have cohomology classes $[\omega]=\alpha$ 
and $[\Omega]=\beta$. The manifolds $X_{\alpha,\beta}$ for 
$(\alpha,\beta)\in U$ form a family diffeomorphic to~$X\times U$.

Next, Kronheimer \cite{Kron2} showed that if $Y$ is any ALE space
asymptotic to $\C^2/H$ then the associated pair $(\alpha,\beta)$
must lie in $U$, and that if $Y_1$ and $Y_2$ are two ALE spaces 
asymptotic to $\C^2/H$ that both yield the same pair $(\alpha,\beta)$, 
then $Y_1,Y_2$ are isomorphic as ALE spaces. Combining these results 
we see that every ALE space $Y$ asymptotic to $\C^2/H$ is isomorphic
to $X_{\alpha,\beta}$ for some pair $(\alpha,\beta)\in U$, so we
have a complete description of all ALE spaces.

The group $\Aut(\Gamma)\ltimes W$ associated to $\C^2/H$ acts on 
Kronheimer's construction, in the following way. The obvious action 
of $\Aut(\Gamma)\ltimes W$ on $H^2(X,\R)\times H^2(X,\C)$ preserves 
the subset $U$, so that $\Aut(\Gamma)\ltimes W$ also acts on $U$. The 
action extends naturally to the hyperk\"ahler quotient construction 
that Kronheimer uses to construct $X_{\alpha,\beta}$, and this shows 
that if $w\in\Aut(\Gamma)\ltimes W$ and $(\alpha,\beta)\in U$, then 
$X_{w\cdot\alpha,w\cdot\beta}$ is isomorphic to $X_{\alpha,\beta}$ 
as an ALE space.

Moreover, if $w\in W$ rather than $\Aut(\Gamma)\ltimes W$, 
then there is a {\it unique} ALE space isomorphism between 
$X_{w\cdot\alpha,w\cdot\beta}$ and $X_{\alpha,\beta}$ that is 
asymptotic to the identity at infinity. One can also show that 
if $X_{\alpha',\beta'}$ is isomorphic to $X_{\alpha,\beta}$ as 
an ALE space, then $(\alpha',\beta')=(w\cdot\alpha,w\cdot\beta)$ 
for some $w\in\Aut(\Gamma)\ltimes W$. If in addition the isomorphism 
between $X_{\alpha',\beta'}$ and $X_{\alpha,\beta}$ is asymptotic 
to the identity at infinity, then~$w\in W$.

Let $w\in\Aut(\Gamma)\ltimes W$. Now $X_{\alpha,\beta}$ and 
$X_{w\cdot\alpha,w\cdot\beta}$ are both diffeomorphic to $X$, under 
diffeomorphisms that are natural up to isotopy. The identification 
between $X_{\alpha,\beta}$ and $X_{w\cdot\alpha,w\cdot\beta}$ that 
comes from their isomorphism as ALE spaces can thus be thought of as 
a diffeomorphism of $X$, up to isotopy. The corresponding isotopy 
class of diffeomorphisms of $X$ is identified with $w\in\Aut(\Gamma)
\ltimes W$, regarding $\Aut(\Gamma)\ltimes W$ as a group of isotopy 
classes of diffeomorphisms of $X$, as above. In particular, the 
identification between $X_{\alpha,\beta}$ and 
$X_{w\cdot\alpha,w\cdot\beta}$ induces the action of $w$ on $H^2(X,\R)$ 
and $H^2(X,\C)$, and this is why it is possible for two isomorphic 
K\"ahler forms apparently to have two different cohomology classes 
$\alpha$ and~$w\cdot\alpha$.

Here is a heuristic description of what is going on. When we 
desingularize $\C^2/H$ we replace the singular point by a
bunch of 2-spheres, and this introduces nontrivial homology
classes in $H_2(X,\Z)$. The Weyl group $W$ then acts as a 
kind of `internal symmetry group' on the new homology classes;
we can visualize elements of $W$ as diffeomorphisms of $X$
that are the identity outside a small neighbourhood of the 
2-spheres. Elements of $\Aut(\Gamma)$ also act as diffeomorphisms
of $X$, but they act nontrivially near infinity.

The ALE spaces $X_{\alpha,\beta}$ for $(\alpha,\beta)\in U$ can 
be thought of as a family of K\"ahler structures upon the fixed real 
4-manifold $X$. Then $X_{\alpha,\beta}$ and $X_{w\cdot\alpha,w\cdot\beta}$ 
represent K\"ahler structures on $X$ that are equivalent under a 
diffeomorphism $\phi$ of $X$ corresponding to $w$. If $w\ne 1$ then 
$\phi$ is not isotopic to the identity, and acts nontrivially on 
$H^2(X,\R)$ and~$H^2(X,\C)$.

\section{Desingularizing Calabi-Yau orbifolds}

We are interested in desingularizing Calabi-Yau orbifolds of 
dimension 3 whose singularities are modelled upon $\C^3/G$, 
where $G$ is a finite subgroup of $SU(3)$, and $\C^3/G$ has 
singularities in codimension two. If we first understand the 
different ways of desingularizing such $\C^3/G$, this will give 
us a local model for how to desingularize more general Calabi-Yau 
orbifolds~$X/G$. 

Suppose $G\subset SU(3)$ is finite and $\C^3/G$ has codimension
two singularities. Pick $x\in\C^3$ such that $xG$ is a generic 
point in the codimension two singular set, and let $H$ be 
$\{h\in G:h(x)=x\}$, the stabilizer subgroup of $x$ in $G$. Then 
there is a natural orthogonal splitting $\C^3=\C\oplus\C^2$,
such that $x\ne 0$ lies in $\C$, and $H$ fixes $\C$ and acts on 
$\C^2$ as a finite, nontrivial subgroup of $SU(2)$. We may write 
$\C^3/H=\C\times(\C^2/H)$, where $\C^2/H$ is one of the {\it 
Kleinian singularities} of~\S 2.

We shall restrict our attention to the case that $H$ {\it is a 
normal subgroup of\/} $G$. If $H$ is not normal then things are
more difficult. So suppose that $H$ is normal in $G$, and let $K$ 
be the quotient group $G/H$. Then $K$ acts naturally on 
$\C\times\C^2/H$, and $\bigl(\C\times\C^2/H\bigr)/K=\C^3/G$. 
Let the notation $X,\Gamma,\Delta,W,\Aut(\Gamma),\rho,U$ and 
$X_{\alpha,\beta}$ all be as defined in the previous section.

We begin by constructing two natural group homomorphisms 
$\phi:K\rightarrow U(1)$ and $\psi:K\rightarrow\Aut(\Gamma)$. Since 
$H$ is the subgroup of $G$ fixing $\C$ and is normal in $G$, it 
follows that $G$ preserves the splitting $\C^3=\C\oplus\C^2$. 
Therefore $G$ is a subgroup of $S(U(1)\times U(2))$, the subgroup 
of $SU(3)$ preserving this splitting, and we may write each element 
$g\in G$ as a pair $(\sigma,\tau)$, where $\sigma\in U(1)$, 
$\tau\in U(2)$, and $\sigma\cdot\det\tau=1$. Then $g\in H$ if and 
only if $\sigma=1\in U(1)$. Define a map $\phi:K\rightarrow U(1)$ by 
$\phi(gH)=\sigma$ for each $g\in G$, where $g=(\sigma,\tau)$. It is 
easy to see that $\phi$ is well-defined, and a group homomorphism.

Write $C_h$ for the conjugacy class of $h$ in $H$, and let $S_H$ be 
the set of nonidentity conjugacy classes in $H$. As $H$ is a normal
subgroup, $g\,C_hg^{-1}$ is also a conjugacy class in $H$ for each
$g\in G$, which is the identity if and only if $C_h$ is. Thus 
$C_h\mapsto gC_hg^{-1}$ defines a map from $S_H$ to itself. Define 
a map from $K\times S_H$ to $S_H$ by $(gH,C_h)\mapsto g\,C_hg^{-1}$.
Then this map is well-defined and is an {\it action} of $K$ on $S_H$.
But there is a natural correspondence between $S_H$, the nonidentity
conjugacy classes in $H$, and the vertices of the Dynkin diagram 
$\Gamma$. Thus $K$ acts on the vertices of $\Gamma$. In fact $K$ acts 
by automorphisms of the whole graph, and this defines the group 
homomorphism $\psi:K\rightarrow\Aut(\Gamma)$ that we want.

Now let $(\alpha,\beta)\in U$, so that $X_{\alpha,\beta}$ is a 
nonsingular ALE space, diffeomorphic to $X$, and asymptotic to 
$\C^2/H$ as in \S 2. Then $\C\times X_{\alpha,\beta}$ desingularizes 
$\C\times\C^2/H$, and has a natural Calabi-Yau structure. Our goal 
is to choose $(\alpha,\beta)$ such that $\C\times X_{\alpha,\beta}$ 
admits a $K$-action preserving this Calabi-Yau structure, which is 
asymptotic to the natural action of $K$ on $\C\times\C^2/H$. To achieve 
this, we must work out what conditions $\alpha$ and $\beta$ must satisfy 
for such a $K$-action to exist.

First consider how $K$ can act on $X$, as a group of diffeomorphisms. 
The action of $K$ on $\C^2/H$ only determines how $K$ should act 
`near infinity' in $X$. But elements of $W$ may be visualized as 
diffeomorphisms of $X$ that are the identity near infinity. Thus 
the action of $K$ on $X$ may not be uniquely determined by the 
asymptotic conditions on it, but instead there may be several such 
actions, differing only by elements of $W$. The data we need to 
determine how $K$ acts on $X$ is a group homomorphism 
$\chi:K\rightarrow\Aut(\Gamma)\ltimes W$, such that $\pi\circ\chi=\psi$, 
where $\pi:\Aut(\Gamma)\ltimes W\rightarrow\Aut(\Gamma)$ is the natural 
projection, so that $\chi$ lifts $\psi$ from $\Aut(\Gamma)$ to 
$\Aut(\Gamma)\ltimes W$. Choose such a homomorphism~$\chi$. 

There always exists at least one such homomorphism,
because there is a canonical choice for $\chi$, which will 
yield desingularizations of $\C^3/G$ with the topology of 
crepant resolutions in the construction below. As $\Aut(\Gamma)$ 
is naturally isomorphic to a subgroup of $\Aut(\Gamma)\ltimes W$,
we can regard $\psi:K\rightarrow\Aut(\Gamma)$ as a homomorphism 
$K\rightarrow\Aut(\Gamma)\ltimes W$, and this is the canonical choice 
for $\chi$. However, for many groups $G,H$ there are other, different 
choices for $\chi$, and we may be able to use these to construct 
resolutions of $\C^3/G$ which do not have the topology of a crepant 
resolution.

Since $\chi:K\rightarrow\Aut(\Gamma)\ltimes W$ is a group homomorphism 
and $\rho$ is a representation of $\Aut(\Gamma)\ltimes W$ on $H^2(X,\R)$
and $H^2(X,\C)$, we see that $\rho\circ\chi$ is a representation
of $K$ on $H^2(X,\R)$ and $H^2(X,\C)$. So, suppose for the moment
that $K$ acts on $X$ as a group of diffeomorphisms, with action 
asymptotic to the prescribed action of $K$ on $\C^2/H$, such that the 
induced action of $K$ on $H^2(X,\R)$ and $H^2(X,\C)$ is~$\rho\circ\chi$.

Next, we choose $(\alpha,\beta)\in U$ and identify $X_{\alpha,\beta}$ 
with $X$ as a real 4-manifold, so that $K$ acts on $X_{\alpha,\beta}$, 
and so on $\C\times X_{\alpha,\beta}$. What is the condition on the pair
$(\alpha,\beta)$ for this $K$-action to preserve the natural Calabi-Yau 
structure on $\C\times X_{\alpha,\beta}$? Let the K\"ahler form and 
holomorphic volume form of $\C$ be $\omega$ and $\Omega$, and let the 
K\"ahler form and holomorphic volume form of $X_{\alpha,\beta}$ be 
$\omega'$ and $\Omega'$, respectively. Then the K\"ahler form of 
$\C\times X_{\alpha,\beta}$ is $\omega+\omega'$, and the holomorphic 
volume form of $\C\times X_{\alpha,\beta}$ is $\Omega\wedge\Omega'$. 
Thus the $K$-action on $\C\times X_{\alpha,\beta}$ must preserve both 
$\omega+\omega'$ and~$\Omega\wedge\Omega'$.

Write $g\in G$ as a pair $(\sigma,\tau)$ as above, where $\sigma\in U(1)$ 
and $\tau\in U(2)$. Then $gH\in K$ acts on $\omega$ and $\Omega$ by 
$gH\cdot\omega=\omega$ and $gH\cdot\Omega=\sigma\cdot\Omega$. Therefore, 
the Calabi-Yau structure of $\C\times X_{\alpha,\beta}$ is $K$-invariant 
if $gH\cdot\omega'=\omega'$ and $\sigma\cdot(gH\cdot\Omega')=\Omega'$. 
Now the cohomology classes of $\omega'$ and $\Omega'$ are $\alpha$ and 
$\beta$ respectively, and $\sigma=\phi(gH)$ from above. Therefore a 
necessary condition on the pair $(\alpha,\beta)$ for the Calabi-Yau 
structure on $\C\times X_{\alpha,\beta}$ to be $K$-invariant is
\begin{equation}
\!\!\!\!\!\!\!\!\!\!\!
\rho\circ\chi(gH)\alpha=\alpha\quad\text{and}\quad
\phi(gH)\cdot\rho\circ\chi(gH)\beta=\beta\quad
\text{for all $gH\in K$.}
\label{abeq}
\end{equation}

It turns out that equation \eq{abeq} is also a {\it sufficient} 
condition for there to exist a $K$-action on $X_{\alpha,\beta}$ with
all the properties we require. In particular, we do not need the
assumption we made above about the existence of a suitable action
of $K$ on $X$ by diffeomorphisms, because equation \eq{abeq}
guarantees this. We now explain why \eq{abeq} is a sufficient 
condition. Recall from \S 2 that if $w\in\Aut(\Gamma)\ltimes W$ 
and $(\alpha,\beta)\in U$ then $X_{\alpha,\beta}$ and 
$X_{w\cdot\alpha,w\cdot\beta}$ are isomorphic as ALE spaces, and that 
if $w\in W$ then there is a unique isomorphism which is asymptotic 
to the identity at infinity.

Extending the arguments used to show this, one can prove the following 
result. Let $gH\in K$ and $(\alpha,\beta)$ and $(\alpha',\beta')$ lie 
in $U$, and consider a map $\Theta_{gH}:X_{\alpha,\beta}\rightarrow 
X_{\alpha',\beta'}$ that is an isomorphism of K\"ahler manifolds, is 
asymptotic to the action of $gH$ on $\C^2/H$, multiplies holomorphic 
volume forms by $\phi(gH)^{-1}$, and acts on cohomology by 
$\rho\circ\chi(gH)$. Then the necessary and sufficient condition for 
there to exist such a map $\Theta_{gH}$ is that $\rho\circ\chi(gH)
\alpha=\alpha'$ and $\phi(gH)\cdot\rho\circ\chi(gH)\beta=\beta'$, and 
if it exists then $\Theta_{gH}$ is unique.

But any solution $(\alpha,\beta)$ to \eq{abeq} satisfies these 
conditions with $(\alpha',\beta')=(\alpha,\beta)$. Therefore, this 
result guarantees the existence and uniqueness of a map $\Theta_{gH}$ 
from $X_{\alpha,\beta}$ to itself with the properties above, for each 
$gH\in K$. It is then easy to show that the maps $\Theta_{gH}$ yield 
an action of $K$ on $X_{\alpha,\beta}$ with all the properties we need. 
We summarize our progress so far in the following Theorem; the final 
part is left as an exercise for the reader.

\begin{thm} Using the above notation, suppose that\/ $\chi:K\rightarrow
\Aut(\Gamma)\ltimes W$ is a group homomorphism such that\/ $\pi\circ\chi
=\psi$, and suppose that\/ $(\alpha,\beta)\in U$ satisfies the condition
that\/ $\rho\circ\chi(gH)\alpha=\alpha$ and\/ $\phi(gH)\cdot\rho\circ\chi
(gH)\beta=\beta$ for all\/ $gH\in K$. Then there exists a unique action
of\/ $K$ on $\C\times X_{\alpha,\beta}$ that preserves the Calabi-Yau
structure of\/ $\C\times X_{\alpha,\beta}$ and is asymptotic to the
natural action of\/ $K$ on $\C\times\C^2/H$. The representation
of\/ $K$ on $H^2(X_{\alpha,\beta},\R)$ induced by this action is
$\rho\circ\chi$, identifying $X_{\alpha,\beta}$ and\/ $X$ as real 
4-manifolds. Moreover, if\/ $Y$ is any ALE space asymptotic to 
$\C^2/H$ such that\/ $\C\times Y$ admits a $K$-action preserving 
the Calabi-Yau structure and asymptotic to the given action on 
$\C\times\C^2/H$, then $Y$ arises from this construction.
\label{ghthm}
\end{thm}

Let us now consider the condition \eq{abeq} more closely. What it 
really means is that $\alpha$ has to be invariant under the action 
$\rho\circ\chi$ of $K$ on $H^2(X,\R)$, but $\beta$ has to be invariant 
under the action $\phi\cdot\rho\circ\chi$ of $K$ on $H^2(X,\C)$. When 
$\phi$ is nontrivial, these two $K$-actions are different, and will 
have different invariant subspaces.

Our construction only works if we are able to choose 
$(\alpha,\beta)\in U$ satisfying \eq{abeq}. Thus by definition of 
$U$, we must find $K$-invariant elements $\alpha$ and $\beta$ such 
that either $\alpha(\delta)\ne 0$ or $\beta(\delta)\ne 0$ for each 
$\delta\in\Delta$. To satisfy this condition in examples, the fact 
that the two $K$-actions are different is important. For instance, 
it can happen that for some $\delta\in\Delta$, every element 
$\alpha\in H^2(X,\R)$ invariant under the $K$-action $\rho\circ\chi$ 
satisfies $\alpha(\delta)=0$. But because $\beta$ must be invariant 
under a different $K$-action $\phi\cdot\rho\circ\chi$, there may 
still exist a suitable element $\beta$ with~$\beta(\delta)\ne 0$.

The goal of this section is to find Calabi-Yau desingularizations of 
quotient singularities $\C^3/G$. We divide the problem into two stages, 
firstly to desingularize $\C^3/H$ to get $\C\times X_{\alpha,\beta}$ 
with a $K$-action, and then secondly to divide by $K$ and desingularize 
the result $\bigl(\C\times X_{\alpha,\beta}\bigr)/K$. So far we have 
discussed only the first stage in this process. But what about the 
second stage?

In fact, at the second stage three things can happen. Firstly, the 
$K$-action on $X_{\alpha,\beta}$ may have no fixed points. In this case 
$\bigl(\C\times X_{\alpha,\beta}\bigr)/K$ has no singularities, and is 
itself a desingularization of $\C^3/G$. Secondly, the singularities
of $\bigl(\C\times X_{\alpha,\beta}\bigr)/K$ may be isolated points.
In this case we must desingularize using a crepant resolution. And 
thirdly, $\bigl(\C\times X_{\alpha,\beta}\bigr)/K$ may have 
singularities in codimension 2. In this case we can of course use a 
crepant resolution, but we are also free to apply the method above again. 

That is, the singularities of $\bigl(\C\times X_{\alpha,\beta}\bigr)/K$ 
are locally modelled on $\C^3/G'$, where $G'$ is a finite subgroup of 
$SU(3)$ that is isomorphic to a subgroup of $K$. As the singularities 
have codimension 2, there is a subgroup $H'$ in $G'$ contained in some 
$SU(2)\subset SU(3)$. We may use the method above to find ways to 
desingularize $\C^3/G'$, and use these as a local model to desingularize 
$\bigl(\C\times X_{\alpha,\beta}\bigr)/K$. Note that $\md{G'}\le\md{K}$ 
as $G'$ is a subgroup of $K$, and $\md{K}=\md{G}/\md{H}<\md{G}$, so that 
$\md{G'}<\md{G}$. Thus, if we use the method iteratively the size of the 
quotient groups decreases at each stage, and the process must terminate.

\section{An example}

We now apply the theory of the previous section to the example
of $\C^3/\Z_4$. First in \S 4.1 we explain how to desingularize 
$\C^3/\Z_4$ in two topologically distinct ways, and then in \S 4.2 
we use this to desingularize a compact Calabi-Yau orbifold~$T^6/\Z_4$.

\subsection{Two ways to desingularize ${\mathbb C}^3/{\mathbb Z}_4$}

Let $\C^3$ have complex coordinates $(z_1,z_2,z_3)$, and define 
$\kappa:\C^3\rightarrow\C^3$ by
\begin{equation}
\kappa:(z_1,z_2,z_3)\longmapsto(-z_1,iz_2,iz_3).
\label{daleq}
\end{equation}
Define $G=\{1,\kappa,\kappa^2,\kappa^3\}$, so that $G$ is a 
finite subgroup of $SU(3)$ isomorphic to $\Z_4$. The only fixed 
point of $\kappa$ and $\kappa^3$ is $(0,0,0)$, but the fixed points 
of $\kappa^2$ are $(z_1,0,0)$ for all $z_1\in\C$. Therefore $\C^3/G$ 
has singularities of codimension two. The subgroup of $G$ fixing the 
points $(z_1,0,0)$ is $H=\{1,\kappa^2\}$, which is a normal subgroup 
of $G$, and preserves the obvious splitting $\C^3=\C\oplus\C^2$. 
Thus the theory of \S 3 applies to $G$ and~$H$.

The Kleinian singularity $\C^2/H$ is $\C^2/\{\pm1\}$. The crepant 
resolution $X$ of $\C^2/\{\pm1\}$ has $H_2(X,\Z)=\Z$. The Dynkin
diagram $\Gamma$ is $A_1$, with $\Aut(\Gamma)=\{1\}$ and Weyl group
$W=\Z_2$. Let $W=\{1,\lambda\}$. Then the generator $\lambda$ of $W$ 
acts on $H_2(X,\Z)$ by multiplication by $-1$. The quotient group 
$K=G/H$ is $\Z_2$ with generator $\kappa H$. Thus the homomorphism 
$\chi:K\rightarrow\Aut(\Gamma)\ltimes W$ of \S 3 maps $\Z_2$ to 
$\Z_2$, and the condition $\pi\circ\chi=\psi$ on $\chi$ is trivial 
since $\Aut(\Gamma)=\{1\}$. Therefore there are two possibilities 
for $\chi$, given by
\begin{equation}
(a)\; \chi(H)=1,\; \chi(\kappa H)=1,\quad\text{and}\quad
(b)\; \chi(H)=1,\; \chi(\kappa H)=\lambda.
\label{chi1eq}
\end{equation}

Since $H^2(X,\R)\cong\R$ and $H^2(X,\C)\cong\C$, the ALE spaces
asymptotic to $\C^2/H$ are parametrized by pairs $(\alpha,\beta)$
with $\alpha\in\R$ and $\beta\in\C$. The condition for 
$(\alpha,\beta)\in U$ is that either $\alpha\ne 0$ or $\beta\ne 0$. 
Let us calculate the conditions \eq{abeq} on $(\alpha,\beta)$ for 
the ALE space $X_{\alpha,\beta}$ to admit a suitable $K$-action, 
for each possibility $(a)$ and $(b)$ in \eq{chi1eq}. These conditions 
involve $\phi$, which is given by $\phi(H)=1$, $\phi(\kappa H)=-1$. 
In case $(a)$, with $gH=\kappa H$, equation \eq{abeq} gives 
$\alpha=\alpha$ and $\beta=-\beta$, which holds if $\beta=0$. In 
case $(b)$ with $gH=\kappa H$, equation \eq{abeq} gives 
$\alpha=-\alpha$ and $\beta=\beta$, which holds if~$\alpha=0$. 

Thus, section 3 gives two different ways $(a)$ and $(b)$ to choose an 
ALE space $X_{\alpha,\beta}$ asymptotic to $\C^2/\{\pm1\}$ together 
with a $K$-action on $\C\times X_{\alpha,\beta}$ asymptotic to the 
given action of $K$ on $\C^3/H$. The next step is to desingularize
$\bigl(\C\times X_{\alpha,\beta}\bigr)/K$ to get a desingularization 
$Y$ of $\C^3/\Z_4$. Here is what happens in each case.

\begin{rlist}
\item[(a)] Let $\alpha\in\R$ be nonzero. Then $K$ acts on $\C\times 
X_{\alpha,0}$. The fixed points of $\kappa H$ in $\C\times X_{\alpha,0}$ 
are a copy of $\CP^1$. Thus $\bigl(\C\times X_{\alpha,0}\bigr)/K$ has 
singularities in codimension two. These singularities admit no 
deformations, but they do have a unique crepant resolution $Y_1$, 
which can be described explicitly using toric geometry. The Betti 
numbers $b^j=b^j(Y_1)$ of $Y_1$ are
\begin{equation}
b^0=b^2=b^4=1,\qquad
b^1=b^3=b^5=b^6=0.
\label{y1betti}
\end{equation}

\item[(b)] Let $\beta\in\C$ be nonzero. Then $X_{0,\beta}$ is
isomorphic as a complex surface to the hypersurface 
$x_1x_3-x_2^2=\beta$ in $\C^3$. Using these coordinates on 
$X_{0,\beta}$, the $K$-action on $\C\times X_{0,\beta}$ is given
by $\kappa H\cdot(z_1,x_1,x_2,x_3)=(-z_1,-x_1,-x_2,-x_3)$. Now 
this action has {\it no fixed points} in $\C\times X_{0,\beta}$, 
since $(0,0,0,0)$ does not satisfy $x_1x_3-x_2^2=\beta$. Thus
$Y_2=\bigl(\C\times X_{0,\beta}\bigr)/K$ is already nonsingular,
and is a desingularization of~$\C^3/\Z_4$.

A careful analysis shows that $Y_2$ retracts onto the subset
\begin{equation}
\bigl\{\pm(0,x_1,x_2,x_3)\in Y_2:
\md{x_1}^2+2\md{x_2}^2+\md{x_3}^2=2\md{\beta}\bigr\},
\end{equation}
which is a copy of $\mathbb{RP}^2$. Thus the fundamental group and 
cohomology of $Y_2$ and $\mathbb{RP}^2$ are isomorphic. So 
$\pi_1(Y_2)\cong\Z_2$, and the Betti numbers $b^j=b^j(Y_2)$ are
\begin{equation}
b^0=1,\qquad b^1=\cdots=b^6=0.
\label{y2betti}
\end{equation}
\end{rlist}

\noindent From \eq{y1betti} and \eq{y2betti} we see that methods 
$(a)$ and $(b)$ yield desingularizations $Y_1$ and $Y_2$ of 
$\C^3/\Z_4$ with rather different topology. The reason for this
difference is that in case $(a)$, $\kappa H$ acts trivially on 
$H_2(X_{\alpha,0},\Z)$, but in case $(b)$, $\kappa H$ acts on 
$H_2(X_{0,\beta},\Z)$ by multiplication by $-1$, so the two 
$K$-actions are topologically distinct.

\subsection{An orbifold $T^6/{\mathbb Z}_4$ and how to desingularize it}

Let $\C^3$ have complex coordinates $(z_1,z_2,z_3)$, and define a 
{\it lattice} $\Lambda$ in $\C^3$ by
\begin{equation}
\Lambda=\bigl\{(a_1+ib_1,a_2+ib_2,a_3+ib_3):a_j,b_j\in\Z\bigr\}.
\end{equation}
Then $\C^3/\Lambda$ is a 6-torus $T^6$, equipped with a flat Calabi-Yau
structure. Let $\kappa$ act on $T^6$ by
\begin{equation}
\kappa:(z_1,z_2,z_3)+\Lambda\longmapsto(-z_1,iz_2,iz_3)+\Lambda,
\label{kaeq}
\end{equation}
as in \eq{daleq}. Then $\kappa$ is well-defined and preserves the
Calabi-Yau structure on $T^6$. Let $G=\{1,\kappa,\kappa^2,\kappa^3\}$ 
be the group generated by $\kappa$, so that $G\cong\Z_4$. Then $T^6/G$ 
is a compact Calabi-Yau orbifold. To understand the singular set 
of $T^6/G$, we shall first find the fixed points of $\kappa,\kappa^2$
and~$\kappa^3$.

The subset of $T^6$ fixed by $\kappa$ and $\kappa^3$ turns out to be 
the 16 points
\begin{equation}
\!\!\!\!\!\!\!\!
\bigl\{(z_1,z_2,z_3)+\Lambda:\;
z_1\in\{0,\ha,\ha i,\ha+\ha i\},\;
z_2,z_3\in\{0,\ha+\ha i\}\bigr\}.
\end{equation}
And the subset of $T^6$ fixed by $\kappa^2$ is 16 copies of $T^2$,
given by
\begin{equation}
\bigl\{(z_1,z_2,z_3)+\Lambda:\; z_1\in\C,\;
z_2,z_3\in\{0,\ha,\ha i,\ha+\ha i\}\bigr\}.
\end{equation}
Twelve of the 16 copies of $T^2$ fixed by $\kappa^2$ are identified
in pairs by the action of $\kappa$, and these contribute 6 copies of
$T^2$ to the singular set of $T^6/G$. On the remaining 4 copies 
$\kappa$ acts as $-1$, so these contribute 4 copies of $T^2/\{\pm1\}$ 
to the singular set. Each $T^2/\{\pm1\}$ contains 4 of the 16 points 
fixed by~$\kappa$.

Therefore the singular set of $T^6/G$ consists of 6 copies of $T^2$, 
with singularities modelled on $T^2\times\C^2/\{\pm1\}$, and 4 copies 
of $T^2/\{\pm1\}$. Each $T^2/\{\pm1\}$ has 4 special points where 
the singularity is modelled on 0 in $\C^3/G$ as in \S 4.1,
and the other singular points look locally like the singularities
of~$\C\times\C^2/\{\pm1\}$. 

To desingularize $T^6/G$, each copy of $T^2$ in the singular set 
can be resolved with a crepant resolution, which replaces the $T^2$ 
by $T^2\times\CP^1$. But each copy of $T^2/\{\pm1\}$ in the singular 
set, may be desingularized using either method $(a)$ or method $(b)$ 
above. For each $k=0,\dots,4$, let $Z_k$ be one of the manifolds 
obtained by desingularizing $T^6/G$ using a crepant resolution 
for each $T^2$ in the singular set, using method $(a)$ for $k$ 
of the $T^2/\{\pm1\}$'s, and using method $(b)$ for the remaining 
$4\!-\!k$ copies of~$T^2/\{\pm1\}$. Then each $Z_k$ is a compact,
nonsingular manifold carrying a family of Calabi-Yau structures.

We shall find the Betti numbers of $Z_k$. The Betti numbers of
$T^6/G$ are
\begin{equation}
\!\!\!\!\!\!\!\!\!\!\!\!
b^0(T^6/\Z_4)=1,\;\> b^1(T^6/\Z_4)=0,\;\>
b^2(T^6/\Z_4)=5,\;\> b^3(T^6/\Z_4)=4.
\label{orbbetti}
\end{equation}
To find the Betti numbers of $Z_k$ we must add on contributions from 
each component of the singular set. The resolution of each copy of $T^2$ 
in the singular set adds 1 to $b^2$ and 2 to $b^3$. Desingularizing a 
$T^2/\{\pm1\}$ using method $(a)$ adds 5 to $b^2$ and fixes $b^3$, but 
desingularizing using method $(b)$ fixes $b^2$ and adds 2 to $b^3$. All 
three processes fix $b^0$ and~$b^1$. 

Thus we calculate that the Calabi-Yau manifolds $Z_k$ have Betti numbers 
\begin{equation}
\!\!\!\!\!\!\!
b^0(Z_k)=1,\;\> b^1(Z_k)=0,\;\> b^2(Z_k)=11+5k,\;\> b^3(Z_k)=24-2k,
\end{equation}
giving Euler characteristic $\chi(Z_k)=12k$. For $k=1,2,3,4$ one can 
show that $Z_k$ is simply-connected, and carries metrics with holonomy 
$SU(3)$. But $Z_0$ is the quotient of $T^2\times K3$ by a free 
$\Z_2$-action, and has fundamental group $\Z_2\ltimes\Z^2$ and 
holonomy~$\Z_2\times SU(2)$.

We have found five compact Calabi-Yau manifolds $Z_0,\dots,Z_4$
that desingularize the same orbifold $T^6/\Z_4$. Now the physicists'
formula \eq{peuler} for the Euler characteristic of desingularizations
of $T^6/\Z_4$ predicts the value 48, which is true when $k=4$ but
false when $k=0,\dots,3$. This is consistent, since $Z_4$ is a crepant
resolution of $T^6/\Z_4$ and it is known that \eq{peuler} holds for 
crepant resolutions, but $Z_0,\dots,Z_3$ are only smoothings of 
$T^6/\Z_4$. Similarly, the Hodge number formulae of Vafa and Zaslow 
hold for $Z_4$, but not for~$Z_0,\dots,Z_3$.

Our examples show that the physicists' Euler characteristic 
\eq{peuler}, and other similar formulae, do not always hold when
a Calabi-Yau orbifold is desingularized by deformation. These are
not the first such examples known, as an example was given by
Vafa and Witten \cite[\S 2]{VaWi}, but so far as the author knows
this phenomenon has not been studied before in a systematic way.
We shall suggest in \S 8 how these examples might be
explained using string theory.

So far we have not justified our claim that each manifold $Z_k$ 
carries a family of Calabi-Yau structures, which converge to the 
singular, flat Calabi-Yau structure on $T^6/G$ in an appropriate 
sense. One way to prove this quite explicitly is to use the analytic 
methods in the author's papers \cite{Joyc2}, \cite{Joyc3} on 
exceptional holonomy; see in particular \cite[Ex.~2, p.~350]{Joyc3}.
Minor modifications to the technique are needed to resolve 
singularities of the kind above, and they will be described in 
the author's book \cite{Joyc4}. The problem can also be approached
using algebraic geometry.

\section{Simultaneous resolution of nodes} 

There is another, already well-understood way in which a 
singular Calabi-Yau 3-fold can have several desingularizations 
with different topology, which involves {\it nodes} in 3-folds.
A {\it node}, or {\it ordinary double point}, is a simple kind 
of singularity of Calabi-Yau 3-folds, modelled on the origin in
$$\bigl\{(z_1,z_2,z_3,z_4)\in\C^4:z_1^2+z_2^2+z_3^2+z_4^2=0\bigr\}.$$
There are two ways to desingularize this:
\begin{rlist}
\item A {\it small resolution} replaces the singular point
with a rational curve $\CP^1$. This can be done in two ways,
related by a {\it flop}.
\item One can {\it deform} away the singularity by changing to
$z_1^2+z_2^2+z_3^2+z_4^2=\epsilon$, for nonzero $\epsilon\in\C$.
Topologically, this replaces the singular point by ${\cal S}^3$.
\end{rlist}

So in case $(i)$ the node is replaced by ${\cal S}^2$, and in 
case $(ii)$ by ${\cal S}^3$. Let $X$ be a singular Calabi-Yau 
manifold with one node. Let $Y_1$ be the real manifold got by 
replacing the node in $M$ by ${\cal S}^2$ as in $(i)$, and let 
$\Sigma_1\subset Y_1$ be the new ${\cal S}^2$. Let $Y_2$ be the 
real manifold got by replacing the node by ${\cal S}^3$ as in 
$(ii)$, and let $\Sigma_2\subset Y_2$ be the new ${\cal S}^3$. 
Then $Y_1$ is Calabi-Yau if and only if $[\Sigma_1]\ne 0$ in 
$H_2(Y_1,\R)$, and $Y_2$ is Calabi-Yau if and only if $[\Sigma_2]\ne 0$ 
in $H_3(Y_2,\R)$. If $X$ is compact and the node the only singular 
point, then exactly one of $[\Sigma_1]$ and $[\Sigma_2]$ is nonzero. 
Thus exactly one of $Y_1$ and $Y_2$ is Calabi-Yau.

Now suppose that $X$ is a singular Calabi-Yau manifold with $k$ 
nodes, for $k>1$. In this case it may be possible to desingularize 
$X$ as a Calabi-Yau manifold in two ways: firstly, by performing a 
small resolution of all the nodes together, and secondly, by 
deforming $X$ so that all the nodes disappear. (This idea was 
originally due to Clemens). For this to happen, the global topology 
of $X$ must satisfy certain conditions. 

Let $Y$ be a small resolution of $X$, let $\Sigma_1,\dots,\Sigma_k
\subset Y$ be the copies of ${\cal S}^2$ introduced at each node, and 
let $[\Sigma_j]$ be the homology classes of the $\Sigma_j$ in $H_2(Y,\R)$. 
Then $Y$ is K\"ahler, and hence a Calabi-Yau 3-fold, if and only if there 
is a class in $H^2(Y,\R)$ that is positive on each $[\Sigma_j]$. Also, 
Friedman \cite[\S 8]{Frie} and Tian \cite{Tian} prove that $X$ admits 
smooth deformations $X_t$ if and only if there exist nonzero constants
$\lambda_1,\dots,\lambda_k$ such that $\lambda_1[\Sigma_1]+\cdots+
\lambda_k[\Sigma_k]=0$ in $H_2(Y,\R)$. This condition is independent 
of the choice of small resolution $Y$, and the deformations $X_t$ are 
Calabi-Yau 3-folds.

Let us now apply these ideas to the case of orbifolds. If $X/G$ 
is a 3-dimensional Calabi-Yau orbifold, there may exist a partial 
desingularization $Y$ of $X/G$, which is nonsingular except for a 
finite number of nodes. One can then try to desingularize $Y$ as
above by small resolutions, deformations, or a combination of both,
to get a number of topologically distinct Calabi-Yau 
desingularizations $Z_1,\dots,Z_k$ of $X/G$, which can have a range
of different Betti numbers. Even in simple examples, the topological 
calculations involved in understanding the possibilities for 
$Z_1,\dots,Z_k$ can be long and difficult.

Here too, the singularities of $X/G$ must be of {\it codimension 
two} for this trick to produce desingularizations $Z_j$ by 
deformation. Let $Y$ be a partial desingularization of $X/G$
with $k$ nodes, let $\tilde Y$ be a small resolution of $Y$,
and let $\Sigma_1,\dots,\Sigma_k$ be the copies of ${\cal S}^2$
in $\tilde Y$ introduced by the resolution. To desingularize
$Y$ by deformation, we need linear relations on the classes
$[\Sigma_j]$ in $H_2(\tilde Y,\R)$. But these relations only 
exist if the singularities of $X/G$ are of codimension two, 
because then the different curves $\Sigma_j$ can be joined 
together by the part of $\tilde Y$ that resolves the codimension 
two singularities. If the singularities of $X/G$ are not of
codimension two, then the nodes in $Y$ are isolated from one 
another, and there are no suitable linear relations on 
the~$[\Sigma_j]$.

We have now proposed two ways in which an orbifold $X/G$
can admit several topologically distinct Calabi-Yau 
desingularizations $Z_1,\dots,Z_k$, firstly the method
of \S 3 using the Weyl group, and secondly the 
method above involving simultaneous resolution of nodes.
What is the relationship between the two? 

In fact the two phenomena are genuinely different, and one 
cannot be explained in terms of the other. The method of
\S 3 is a two stage process, where in the first
stage one desingularizes $\C^3/H$ in a $K$-invariant way,
and we are free to choose some topological data $\chi$.
The method of this section comes in at the second stage, 
since it may give us several ways to desingularize $\bigl(\C\times 
X_{\alpha,\beta}\bigr)/K$. In particular, the data $\chi$ still
makes sense on a manifold with nodes, and is not changed
by the choice of how to resolve the nodes.

\section{Another example}

In this section we consider an example that combines the ideas
of \S 3 and \S 5, and shows how complicated the
business of desingularizing orbifolds can be. First in 
\S 6.1 we describe the different ways to desingularize 
$\C^3/\Z_2^2$ as a Calabi-Yau manifold. Then in \S 6.2
we apply this to study the possible desingularizations of
the orbifold $T^6/\Z_2^2$. This turns out to be a complex
problem, which we do not solve completely.

\subsection{Desingularizations of $\C^3/\Z_2^2$}

Let $\C^3$ have complex coordinates $(z_1,z_2,z_3)$, and define 
$\kappa_j:\C^3\rightarrow\C^3$ by
\begin{equation}
\begin{split}
\!\!\!\!\!\!\!\!\!\!\!
\kappa_1:(z_1,z_2,z_3)\mapsto(z_1,-z_2,-z_3),&\quad
\kappa_2:(z_1,z_2,z_3)\mapsto(-z_1,z_2,-z_3)\\
\text{and}&\quad\kappa_3:(z_1,z_2,z_3)\mapsto(-z_1,-z_2,z_3).\\
\end{split}
\label{ka123eq}
\end{equation}
Then $G=\{1,\kappa_1,\kappa_2,\kappa_3\}$ is a subgroup of $SU(3)$ 
isomorphic to $\Z_2^2$. The singular set of $\C^3/G$ splits 
into 3 pieces: the points $\pm(z_1,0,0)$ coming from the
fixed points of $\kappa_1$, the points $\pm(0,z_2,0)$ from the
fixed points of $\kappa_2$, and the points $\pm(0,0,z_3)$ from 
the fixed points of $\kappa_3$. Each piece is a copy of 
$\C/\{\pm1\}$, and they meet at~$(0,0,0)$.

Define $H_j=\{1,\kappa_j\}$ for $j=1,2,3$. Then $H_1,H_2$ and $H_3$ 
are normal subgroups of $G$ with $\C^3/H_j\cong\C\times\C^2/\{\pm1\}$.
The quotient groups $K_j=G/H_j$ are isomorphic to $\Z_2$, and act 
upon $\C^3/H_j$. Thus we can apply the method of \S 3 to 
$\C^3/G$ in 3 different ways, by starting with $H_1,H_2$ or $H_3$, 
and these 3 ways correspond to the 3 pieces of the singular set. 
Now $\C^2/\{\pm1\}$ has Dynkin diagram $\Gamma=A_1$, with $\Aut(\Gamma)
=\{1\}$ and Weyl group $W=\{1,\lambda\}$ isomorphic to~$\Z_2$.

Thus, by \S 3, every Calabi-Yau desingularization $Y$ of $\C^3/\Z_2^2$ 
has three pieces of topological data, the group homomorphisms 
$\chi_j:K_j\rightarrow\Aut(\Gamma)\ltimes W$ for $j=1,2,3$. Here 
$K_1=\{H_1,\kappa_2H_1\}$ and $\Aut(\Gamma)\ltimes W=\{1,\lambda\}$ are 
both isomorphic to $\Z_2$, so there are two possibilities for~$\chi_1$,
\begin{equation}
\!\!\!\!\!\!\!\!\!\!
(a)\;\chi_1(H_1)=1,\;\chi_1(\kappa_2H_1)=1,\text{\ and\ }
(b)\;\chi_1(H_1)=1,\;\chi_1(\kappa_2H_1)=\lambda.
\end{equation}
As a shorthand, we shall write $\chi_1=1$ to denote case $(a)$
and $\chi_1=-1$ to denote case $(b)$. Similarly, there are
two possibilities for each of $\chi_2$ and $\chi_3$, which we
will also write $\chi_2=\pm1$, $\chi_3=\pm1$. We can think of
$\chi_1,\chi_2$ and $\chi_3$ as describing the topology of $Y$
near infinity.

In this section we will study all the different ways to 
desingularize $\C^3/\Z_2^2$ as a Calabi-Yau manifold. First we 
describe the deformations of $\C^3/\Z_2^2$. Let $\gamma:\C^3/\Z_2^2
\rightarrow\C^4$ be given by $\gamma\bigl((z_1,z_2,z_3)G\bigr)= 
(z_1^2,z_2^2,z_3^2,z_1z_2z_3)$. Then $\gamma$ is well-defined, 
and induces an isomorphism between $\C^3/G$ and the hypersurface
\begin{equation}
W_{0,0,0,0}=\bigl\{(x_1,x_2,x_3,x_4)\in\C^4:x_1x_2x_3-x_4^2=0\bigr\}
\end{equation}
in $\C^4$. Let $\alpha,\beta_1,\beta_2$ and $\beta_3$ be complex numbers, 
and define
\begin{equation}
\!\!\!\!\!\!\!\!
W_{\alpha,\beta_1,\beta_2,\beta_3}\!=\!\bigl\{
(x_1,x_2,x_3,x_4)\!\in\!\C^4\!:\!x_1x_2x_3\!-\!x_4^2
\!=\!\alpha\!+\!\beta_1x_1\!+\!\beta_2x_2\!+\!\beta_3x_3\bigr\}.
\end{equation}
Then $W_{\alpha,\beta_1,\beta_2,\beta_3}$ is a {\it deformation} of 
$\C^3/\Z_2^2$. For generic $\alpha,\dots,\beta_3$ the hypersurface 
$W_{\alpha,\beta_1,\beta_2,\beta_3}$ is nonsingular, but for some special 
values of $\alpha,\dots,\beta_3$ it has singularities. If 
$W_{\alpha,\beta_1,\beta_2,\beta_3}$ is singular, then it can be resolved
with a crepant resolution to make it nonsingular. 

Thus, each set of values of $\alpha,\dots,\beta_3$ may give one
or more ways to desingularize $\C^3/\Z_2^2$ as a Calabi-Yau 
manifold. We will now list the different cases that arise in
this way. In each case we will give the values of $\chi_1,\chi_2$
and $\chi_3$, and the Betti numbers $b^2$ and $b^3$ of the 
desingularization.

\begin{rlist}
\item $W_{0,0,0,0}$ is isomorphic to $\C^3/\Z_2^2$. It has 4 
possible crepant resolutions, which are easily described using
toric geometry. Each has $\chi_1=\chi_2=\chi_3=1$ and Betti 
numbers $b^2=3$ and~$b^3=0$.

\item $W_{\alpha,0,0,0}$ for $\alpha\ne 0$. This is nonsingular,
has $\chi_1=\chi_2=\chi_3=1$, and Betti numbers
$b^2=2$ and~$b^3=1$.

\item $W_{0,\beta_1,0,0}$ for $\beta_1\ne 0$. This is isomorphic 
to $\bigl(\C\times X_{0,\beta_1}\bigr)/K_1$, and has singularities 
at the points $(0,x_2,x_3,0)$ for $x_2x_3=\beta_1$. It has a
unique crepant resolution, by blowing up the singular set,
which has $\chi_1=-1$ and $\chi_2=\chi_3=1$, and Betti numbers 
$b^2=1$ and~$b^3=1$.

\item $W_{\alpha,\beta_1,0,0}$ for $\alpha,\beta_1\ne 0$. This is 
nonsingular, and is a smooth deformation of the resolution in $(iii)$, 
with the same topology and values of $\chi_j$ and~$b^k$.

\item[$(v),(vi)$] As $(iii)$ and $(iv)$ but with $\beta_2$
nonzero instead of $\beta_1$, and $\chi_2=-1$ instead of~$\chi_1$.

\item[$(vii),(viii)$] As $(iii)$ and $(iv)$ but with $\beta_3$
nonzero instead of $\beta_1$, and $\chi_3=-1$ instead of~$\chi_1$.

\item[$(ix)$] $W_{\alpha,\beta_1,\beta_2,\beta_3}$ with 
$\beta_1,\beta_2,\beta_3\ne 0$ and $\alpha^2\ne 4\beta_1\beta_2\beta_3$. 
This is nonsingular and has $\chi_1=\chi_2=\chi_3=-1$ and Betti numbers 
$b^2=0$ and~$b^3=1$.
\end{rlist}

For a few special values of $\alpha,\dots,\beta_3$, we cannot resolve
$W_{\alpha,\beta_1,\beta_2,\beta_3}$ as a Calabi-Yau manifold with the
appropriate asymptotic behaviour, so it does not appear on the 
above list. Here are the missing cases, with the reason why.

\begin{itemize}
\item If exactly one of $\beta_1,\beta_2,\beta_3$ is zero, say $\beta_1$, 
then $W_{\alpha,0,\beta_2,\beta_3}$ is nonsingular and is topologically
equivalent to case $(ix)$. However, we should regard it as being 
`singular at infinity'.
\item Also, $W_{\alpha,\beta_1,\beta_2,\beta_3}$ with 
$\beta_1,\beta_2,\beta_3\ne 0$ and $\alpha^2=4\beta_1\beta_2\beta_3$ 
has a single node at $x_j=-\alpha/2\beta_j$ for $j=1,2,3$. However, 
neither of the small resolutions of it are K\"ahler manifolds. 
\end{itemize}

Here is what we mean by `singular at infinity'. Our goal is to
construct Calabi-Yau manifolds that desingularize $\C^3/\Z_2^2$.
As with the ALE spaces of \S 2, we expect these manifolds to be 
asymptotic to $\C^3/\Z_2^2$ at infinity, and the metrics on
them to be asymptotic at infinity to the Euclidean metric on 
$\C^3/\Z_2^2$, in some suitable sense. However, because the
singularities of $\C^3/\Z_2^2$ extend to infinity, things are
more complicated than they seem at first.

I have studied this problem, and I have found a good definition
for the idea of ALE space in the case of non-isolated singularities,
and have also proved an existence result for Calabi-Yau metrics
satisfying this definition, by adapting Yau's proof of the Calabi
conjecture. I hope to publish these results in my forthcoming book
\cite{Joyc4}. In cases $(i)$-$(ix)$ above, my results do guarantee
the existence of Calabi-Yau metrics on the given desingularizations,
for suitable choices of the K\"ahler class.

Roughly speaking, in this case the asymptotic conditions on the 
metrics are as follows. If at least two of $z_1,z_2,z_3$ are very
large, then the metric on the desingularization near the point 
$(z_1,z_2,z_3)G$ in $\C^3/\Z_2^2$ must be close to the flat metric
on $\C^3/\Z_2^2$. But if only one of $z_1,z_2,z_3$ is large, say 
$z_1$, then the metric in the desingularization near the point 
$(z_1,z_2,z_3)G$ in $\C^3/\Z_2^2$ must be close to the product 
Calabi-Yau metric on $\C\times X_{\delta,\epsilon}$, where $z_1$ is the 
coordinate in $\C$, and $X_{\delta,\epsilon}$ is an ALE space asymptotic 
to $\C^2/\{\pm1\}$, which has coordinates~$\pm(z_2,z_3)$.

We say the Calabi-Yau metric is {\it singular at infinity} if
the ALE space $X_{\delta,\epsilon}$ appearing in this asymptotic 
condition is singular -- in this case, if $X_{\delta,\epsilon}=
\C^2/\{\pm1\}$. We have excluded cases like $W_{\alpha,0,\beta_2,
\beta_3}$ for $\beta_2,\beta_3\ne 0$ from our list because they are 
singular at infinity, so that the singularities $\pm(z_1,0,0)$ for 
$z_1$ very large, effectively remain unresolved. It can be shown, 
although we will not prove this, that every desingularization of 
$\C^3/\Z_2^2$ as a Calabi-Yau manifold, that is not `singular at 
infinity', is modelled on one of cases $(i)$-$(ix)$ above. 

In cases $(ii)$, $(iv)$, $(vi)$ and $(viii)$ above, there are 
nontrivial conditions upon the K\"ahler class for the metrics to be 
nonsingular at infinity. The allowed values for the K\"ahler class 
split into several connected components -- six components in case 
$(ii)$ and two components in cases $(iv)$, $(vi)$ and $(viii)$.
The connected component of the K\"ahler class can be regarded as an
extra topological choice in the desingularization; but we will
not discuss this issue here.

The relationship between cases $(i)$ and $(ii)$ may be understood
in terms of the ideas of \S 5. There exists a partial 
resolution (in fact, three different partial resolutions) of
$\C^3/\Z_2^2$ with a single node. This partial resolution can
be resolved by a small resolution, giving one of the four
manifolds in case $(i)$. Alternatively, the node can be
deformed away, giving case~$(ii)$. 

Observe that in the possible desingularizations of $\C^3/\Z_2^2$,
we can have 0,1 or 3 of $\chi_1,\chi_2$ and $\chi_3$ equal to $-1$,
but we cannot have exactly 2 of $\chi_1,\chi_2$ and $\chi_3$ equal
to $-1$. Thus we cannot choose $\chi_1,\chi_2$ and $\chi_3$ 
independently. The moral is that when we desingularize $\C^3/\Z_2^2$
or other orbifolds in which the codimension 2 singularities split 
into several pieces, the topological choices for different pieces 
of the singular set are not in general independent, but are subject
to constraints involving all the pieces.

\subsection{Classifying the desingularizations of $T^6/\Z_2^2$}

Let $\Lambda$ be as in \S 4.2, so that $\C^3/\Lambda$ is a 6-torus 
$T^6$ with a flat Calabi-Yau structure. Let $\kappa_1,\kappa_2$ and 
$\kappa_3$ act on $T^6$ by
\begin{equation}
\begin{split}
\kappa_1&:(z_1,z_2,z_3)+\Lambda\longmapsto(z_1,-z_2,-z_3)+\Lambda,\\
\kappa_2&:(z_1,z_2,z_3)+\Lambda\longmapsto(-z_1,z_2,-z_3)+\Lambda,
\end{split}
\end{equation}
and $\kappa_3=\kappa_1\kappa_2$, as in \eq{ka123eq}. Then 
$G=\{1,\kappa_1,\kappa_2,\kappa_3\}$ acts on $T^6$ preserving the 
Calabi-Yau structure, and is a group isomorphic to $\Z_2^2$. The 
quotient $T^6/G$ is a Calabi-Yau orbifold, with singularities modelled 
on $\C^3/\Z_2^2$. This orbifold was studied by Vafa and Witten 
\cite[\S 2]{VaWi} (see also \cite[\S 4.2]{AsMo}), who showed that 
$T^6/\Z_2^2$ can be resolved by crepant resolution, in many different 
ways, to get a Calabi-Yau manifold with $h^{1,1}=51$ and $h^{2,1}=3$. 
But they also found one way to desingularize $T^6/\Z_2^2$ by deformation, 
to get another Calabi-Yau manifold with $h^{1,1}=3$ and~$h^{2,1}=115$.

Let us now consider how to describe {\it all} the different 
possible ways to desingularize $T^6/\Z_2^2$ as a Calabi-Yau
manifold. We will not be able to offer a complete classification, 
because the calculations involved are extremely complex. However, 
we can explain the first steps in this classification, and we will 
see that there are in fact a large number of different ways to 
desingularize $T^6/\Z_2^2$, of which the possibilities found by 
Vafa and Witten represent two extremes.

We shall regard $T^6$ as a product $T^2\times T^2\times T^2$. Then 
$\Z_2$ acts on each copy of $T^2$, so that $\Z_2^3$ acts on $T^6$, 
and $G$ is a subgroup of this $\Z_2^3$. This $\Z_2$-action on $T^2$ 
has 4 fixed points $p_1,\dots,p_4$. The fixed points of $\kappa_1$ on 
$T^6$ are $T^2\times p_j\times p_k$ for $j,k=1,\dots,4$, which is 16 copies 
of $T^2$. For $j,k=1,\dots,4$, define $A_{jk}=T^2/\Z_2\times p_j\times p_k
\subset T^6/\Z_2^2$. Similarly, for $i,k=1,\dots,4$, define $B_{ik}
=p_i\times T^2/\Z_2\times p_k\subset T^6/\Z_2^2$, and for $i,j=1,\dots,4$ 
define~$C_{ij}=p_i\times p_j\times T^2/\Z_2\subset T^6/\Z_2^2$.

The singular set of $T^6/\Z_2^2$ is the union of these sets $A_{jk}$,
$B_{ik}$ and $C_{ij}$. The $A_{jk}$ come from the fixed points of
$\kappa_1$, the $B_{ik}$ from $\kappa_2$, and the $C_{ij}$ from 
$\kappa_3$. Each of the $A_{jk},B_{ik}$ and $C_{ij}$ is a copy of 
$T^2/\Z_2$. They are not disjoint, but for each $i,j,k=1,\dots,4$ 
the three sets $A_{jk},B_{ik}$ and $C_{ij}$ intersect in the point 
$p_{ijk}=p_i\times p_j\times p_k$ in $T^6/\Z_2^2$. The $p_{ijk}$ are 
the 64 singular points in $T^6/\Z_2^2$ which have a singularity 
modelled on $0$ in~$\C^3/\Z_2^2$.

Now, following the method of \S 3, to desingularize
$T^6/\Z_2^2$ we must first choose some topological data about
the desingularization, the group homomorphism $\chi$. As we saw
above, for the $\C^3/\Z_2^2$ singularity there are 3 pieces of 
data $\chi_1,\chi_2,\chi_3$ corresponding to the $\kappa_1,\kappa_2$
and $\kappa_3$ singularities, and each $\chi_j$ can take the values 
$\pm1$. In our case, a little thought shows that $\chi_1$ gives
topological information about the way the singularities $A_{jk}$
are resolved. One can show that $\chi_1$ must be constant on
each $A_{jk}$, since $A_{jk}$ is connected, but different 
$A_{jk}$ can have different values of $\chi_1$. Write $\chi_{1,jk}$
for the value of $\chi_1$ on $A_{jk}$. Then for $j,k=1,\dots,4$,
we have~$\chi_{1,jk}=\pm1$.

Similarly, $\chi_2$ gives information on how the $B_{ik}$ are 
resolved, and we write $\chi_{2,ik}$ for the value of $\chi_2$ 
on $B_{ik}$, and $\chi_3$ gives information on how the $C_{ij}$ 
are resolved, and we write $\chi_{3,ij}$ for the value of $\chi_3$ 
on $C_{ij}$. Thus, to desingularize $T^6/\Z_2^2$ we must first 
choose the values of $\chi_{1,jk},\chi_{2,ik}$ and $\chi_{3,ij}$. 
These are 48 variables taking the values $\pm 1$, so there are 
$2^{48}$, or about $2\!\cdot\!8\times 10^{14}$ possible choices.

Now, we saw above that the possible desingularizations $(i)$-$(ix)$ 
of $\C^3/\Z_2^2$ as a Calabi-Yau manifold allow $0,1$ or 3 of 
$\chi_1,\chi_2,\chi_3$ to be $-1$, but not two to be $-1$. This 
condition applies at each of the 64 points $p_{ijk}$. Therefore, 
a necessary condition for the data $\chi_{i,jk}$ to represent a 
possible Calabi-Yau desingularization is that for each set of 
values $i,j,k=1,\dots,4$, exactly 0,1 or 3 of $\chi_{1,jk},
\chi_{2,ik}$ and $\chi_{3,ij}$ are $-1$, but not two of them. 

This condition excludes nearly all of the $2^{48}$ choices for 
the $\chi_{i,jk}$, but there are still many choices for which 
this condition is satisfied, although we have not been able to 
count them. However, we can give four explicit families of 
solutions to the conditions, and so find a lower limit for their 
number. For the first family, let $\delta_i,\epsilon_j$ and $\zeta_k$ 
take the values $\pm 1$ for $i,j,k=1,\dots,4$, and define
\begin{equation}
\chi_{1,jk}=\epsilon_j\zeta_k,\qquad\chi_{2,ik}=\delta_i\zeta_k
\qquad\text{and}\qquad\chi_{3,ij}=-\delta_i\epsilon_j.
\label{dezeq}
\end{equation}
Then $\chi_{1,jk}\chi_{2,ik}\chi_{3,ij}=-1$ for all $i,j,k$,
and this means that either 1 or 3 of $\chi_{1,jk},\chi_{2,ik}$ 
and $\chi_{3,ij}$ are equal to $-1$, but not 0 or 2. Conversely,
any set of values of $\chi_{1,jk},\chi_{2,ik}$ and $\chi_{3,ij}$
for which this holds may be written in the form \eq{dezeq}.
There are $2^{12}$ possible values for $\delta_i,\epsilon_j$ and 
$\zeta_k$, but reversing the sign of all the $\delta_i,\epsilon_j$ 
and $\zeta_k$ does not change the $\chi_{i,jk}$, so this gives
$2^{11}=2048$ different solutions for the~$\chi_{i,jk}$.

For the second family, let $\chi_{2,ik}=\chi_{3,ij}=1$ for all
$i,j,k$, and let $\chi_{1,jk}$ be $\pm 1$. Clearly, for all $i,j,k$
either 0 or 1 of $\chi_{1,jk},\chi_{2,ik}$ and $\chi_{3,ij}$ are 
equal to $-1$, so the condition is satisfied. There are 
$2^{16}=65536$ possible choices for the $\chi_{1,jk}$. Similarly,
by putting $\chi_{1,jk}=\chi_{3,ij}=1$, and by putting $\chi_{1,jk}
=\chi_{2,ik}=1$, we get two other families of $2^{16}$ choices.
In total, allowing for repeated choices, we have found 198651
different sets of values for $\chi_{1,jk},\chi_{2,ik}$ and 
$\chi_{3,ij}$ in which the conditions are satisfied. This is a
lower limit on the number of solutions, which is probably rather
larger than this.

Next, having chosen a set of suitable values of the $\chi_{i,jk}$,
we must look for Calabi-Yau desingularizations of $T^6/\Z_2^2$ with 
this data. There are still further topological choices to make. 
At each of the 64 points $p_{ijk}$, we must choose one of the 
desingularizations $(i)$-$(ix)$ above that is consistent with the 
values of $\chi_{1,jk},\chi_{2,ik}$ and $\chi_{3,ij}$ already chosen. 
For instance, if $\chi_{1,jk}=\chi_{2,ik}=\chi_{3,ij}=1$ then either 
case $(i)$ or case $(ii)$ will do. There are also more subtle 
topological choices to do with Weyl groups and the connected 
component of the K\"ahler class, which we will not go into.

Having made all these topological choices, we can finally construct 
a unique real 6-manifold $Y$ that desingularizes $T^6/\Z_2^2$, which 
{\it locally} has the topology of a Calabi-Yau desingularization.
However, many of these 6-manifolds do not admit Calabi-Yau structures
desingularizing $T^6/\Z_2^2$. Recall that in \S 5 we discussed the 
global topological issues involved in desingularizing a Calabi-Yau 
manifold with finitely many nodes. In this situation there are some 
rather similar conditions that must be satisfied for a Calabi-Yau 
structure to exist on~$Y$.

But in some special cases we can see quite easily that the 
Calabi-Yau structures exist. For instance, in the second family
above with $\chi_{2,ik}=\chi_{3,ij}=1$, if we choose 
desingularization $(i)$ for $p_{ijk}$ when $\chi_{1,jk}=1$ and
desingularization $(iii)$ for $p_{ijk}$ when $\chi_{1,jk}=-1$,
then one can prove that the resulting manifold has a Calabi-Yau 
structure, which is a crepant resolution of $(T^2\times K3)/\Z_2$.
Using the same trick with the third and fourth families gives
a total of 196606 different sets of values of the $\chi_{i,jk}$ 
which do correspond to Calabi-Yau desingularizations; and as 
there are four topological choices for resolution $(i)$, these 
will lead to many more manifolds.

Our discussion has shown that the problem of classifying all the
possible Calabi-Yau desingularizations of $T^6/\Z_2^2$ is of great 
complexity. There are a large number of choices to be made, but
these choices are subject to many complicated conditions. The 
author's feeling is that these conditions are not too restrictive, 
and the number of different ways of desingularizing $T^6/\Z_2^2$ 
is probably very large. But $T^6/\Z_2^2$ is in fact one of the 
simplest and most obvious Calabi-Yau orbifolds that one can 
think of, and the problems involved in analyzing more complex
examples must be even worse!

Finally, we explain how the desingularizations of Vafa and Witten
\cite[\S 2]{VaWi} fit into this framework. Crepant resolutions
of $T^6/\Z_2^2$ have all $\chi_{i,jk}=1$, and each point $p_{ijk}$
is resolved using case $(i)$ above. The nonsingular deformation of 
$T^6/\Z_2^2$ given by Vafa and Witten has all $\chi_{i,jk}=-1$, and
each point $p_{ijk}$ is resolved using case $(ix)$ above; this leads
to a unique manifold. These are two extremes in the possible choices
for the $\chi_{i,jk}$, and there are many possibilities in between.

We can also relate another part of Vafa and Witten's ideas to our
analysis above. Vafa and Witten make their desingularization in
two stages, by first deforming to a singular Calabi-Yau manifold
with 64 nodes which they say is a `mirror partner' to the crepant
resolutions of $T^6/\Z_2^2$, and then by deforming away the nodes
to get a nonsingular manifold. Above we explained that if
$\beta_1,\beta_2,\beta_3\ne 0$ and $\alpha^2=4\beta_1\beta_2\beta_3$, 
then $W_{\alpha,\beta_1,\beta_2,\beta_3}$ has a single node, which 
vanishes under deformation. It seems clear that Vafa and Witten's 
singular Calabi-Yau manifold is modelled on this.

\section{Exceptional holonomy}

Calabi-Yau manifolds can be described as Riemannian manifolds
with {\it holonomy group} $SU(3)$. Now the {\it exceptional 
holonomy groups} are two special cases in the classification
of Riemannian holonomy groups, the holonomy groups $G_2$ in 7 
dimensions, and $Spin(7)$ in 8 dimensions. They share many 
properties with the holonomy groups $SU(3)$. In particular,
compact Riemannian 7- and 8-manifolds with holonomy $G_2$ or 
$Spin(7)$ can be made by desingularizing orbifolds $T^7/G$ 
or $T^8/G$ in a special way, and this method was used by the 
author \cite{Joyc1,Joyc2,Joyc3,Joyc4} to construct the first 
known examples. 

The exceptional holonomy groups are also important in String 
Theory -- see for instance Shatashvili and Vafa \cite{ShVa} -- 
in the same way that Calabi-Yau manifolds are. Suppose that 
$X$ is a nonsingular 7- or 8-manifold with holonomy $G_2$ or 
$Spin(7)$, constructed by desingularizing an orbifold $Y/G$.
Then, just as with the Euler characteristic formula \eq{peuler}
and the Hodge number formulae of Vafa and Zaslow for the
case of Calabi-Yau manifolds, String Theory can be used to
predict topological information about $X$ from the orbifold $Y/G$. 
This appears implicitly in Shatashvili and Vafa~\cite{ShVa}.

For a 7-manifold with holonomy $G_2$, the important Betti numbers
are $b^2$ and $b^3$, and String Theory can be used to predict their
sum $b^2+b^3$. For an 8-manifold with holonomy $Spin(7)$, the 
important Betti numbers are $b^2,b^3$ and $b^4=b^4_++b^4_-$, and 
String Theory predicts the linear combinations $2b^2+b^4$, $b^3$ 
and $b^4_+-b^4_-$. In particular, the Euler characteristic formula
\eq{peuler} should apply without change to 8-manifolds with holonomy
$Spin(7)$. (For 7-manifolds the Euler characteristic is of course
zero.)

Now, using the techniques of \cite{Joyc1}-\cite{Joyc4} and the
ideas of \S 3-\S 6, one can construct examples
of orbifolds which have a number of topologically distinct
resolutions with holonomy $G_2$ or $Spin(7)$, for some of which 
the String Theory formulae do not hold. To actually prove this,
one needs to use more sophisticated techniques than those of
\cite{Joyc1}-\cite{Joyc3}, which will be explained in~\cite{Joyc4}.

In the case of Calabi-Yau manifolds, we have a clear distinction
between crepant resolutions, for which the String Theory formulae
are always true (at least in dimensions 2 and 3), and deformations,
for which (in dimensions 3 and above) the formulae are often false. 
However, in the geometry of $G_2$ and $Spin(7)$ the distinction 
between crepant resolutions and deformations no longer makes sense.
So it seems that we cannot separate out a special class of resolutions 
for which the String Theory formulae can be conjectured, or proved,
to hold. Perhaps further developments in String Theory will make
the matter clearer.

We shall now present an example of an {\it isolated} quotient 
singularity $\R^8/G$, which has several resolutions within holonomy 
$Spin(7)$ for which \eq{peuler} does not hold. Identify $\R^8$ with 
$\C^4$, which has complex coordinates $(z_1,z_2,z_3,z_4)$, and the 
standard Euclidean metric and holomorphic volume form. This induces 
a flat $SU(4)$-structure on $\R^8$. But the holonomy groups $SU(4)$ 
and $Spin(7)$ are subgroups of $O(8)$, such that $SU(4)\subset 
Spin(7)\subset O(8)$. Because of this, an $SU(4)$-structure on an
8-manifold induces a unique $Spin(7)$-structure on the same manifold.

Define maps $\kappa,\lambda:\R^8\rightarrow\R^8$ by
\begin{equation}
\!\!\!\!\!\!\!\!
\kappa:(z_1,\dots,z_4)\mapsto(iz_1,iz_2,iz_3,iz_4),
\;\>\lambda:(z_1,\dots,z_4)\mapsto
(\overline z_2,-\overline z_1,\overline z_4,-\overline z_3).
\label{kalaeq}
\end{equation}
Then $\kappa$ and $\lambda$ satisfy the relations $\kappa^4=\lambda^4=1$, 
$\kappa^2=\lambda^2$ and $\kappa\lambda=\lambda\kappa^3$, and they 
generate a group $G$ of automorphisms of $\R^8$, which is nonabelian 
and of order 8. The quotient $\R^8/G$ has one singular point, the origin. 
It can be shown that $\kappa$ lies in both $SU(4)$ and $Spin(7)$, and 
$\lambda$ lies in $Spin(7)$ but not in $SU(4)$. Thus $G$ is a subgroup 
of $Spin(7)$, and the subgroup $H_1=\{1,\kappa,\kappa^2,\kappa^3\}$ of 
$G$ is also a subgroup of $SU(4)$. Note that $H_1$ is a normal subgroup 
of $G$, and $K_1=G/H_1$ is isomorphic to~$\Z_2$.

Following the ideas of \S 3, we will first desingularize $\R^8/H_1$ 
in a $K_1$-invariant way. As $H_1$ is a subgroup of $SU(4)$, we shall 
desingularize $\C^4/H_1$ with holonomy $SU(4)$, that is, as a 
Calabi-Yau manifold. Now $\C^4/H_1$ admits a unique crepant resolution 
$X$, by blowing up the singular point, in which the singular point 
is replaced by a copy of $\CP^3$. But the singularities of $\C^4/H_1$ 
are rigid under deformation, so this crepant resolution is the only 
way to desingularize $\C^4/H_1$ as a Calabi-Yau manifold.

There exists an explicit 1-parameter family of ALE Calabi-Yau 
metrics on $X$, which was written down by Calabi \cite[p.~285]{Cala}. 
Examining these metrics, one can easily verify that the map $\lambda$ of 
\eq{kalaeq} extends to an isometric involution $\lambda':X\rightarrow X$ 
that is asymptotic to $\lambda$ at infinity. This involution has no fixed 
points, so $Y_1=X/\{1,\lambda'\}$ is nonsingular. The Calabi-Yau metric 
on $X$ gives an $SU(4)$-structure, which induces a $Spin(7)$-structure on 
$X$, and $\lambda'$ preserves this $Spin(7)$-structure (but not the 
$SU(4)$-structure). 

Thus, $Y_1$ is a nonsingular 8-manifold asymptotic to $\R^8/G$, and it 
carries a 1-parameter family of ALE metrics and $Spin(7)$-structures
that converge to the orbifold metric on $\R^8/G$. The holonomy group
of these metrics is $\Z_2\ltimes SU(4)$, which is a subgroup of $Spin(7)$.
The fundamental group $\pi_1(Y_1)$ is $\Z_2$, and the Betti numbers
of $Y_1$ are $b^1=b^2=b^3=b^4_+=0$ and $b^4=b^4_-=1$. In particular,
this means that the resolution of the singularity adds 1 to the Euler 
characteristic. But a na\"\i ve application of String Theory ideas
suggests that the resolution should add 5 to the Euler characteristic, 
which is the number of nonidentity conjugacy classes in~$G$.

In fact, we can resolve $\R^8/G$ in three slightly different ways. 
Define $H_2=\{1,\kappa\lambda,\kappa^2,\kappa^3\lambda\}$ and 
$H_3=\{1,\lambda,\lambda^2,\lambda^3\}$. Then $H_2$ and $H_3$ are 
also normal subgroups of $G$, and $K_2=G/H_2$ and $K_3=G/H_3$ are 
isomorphic to $\Z_2$. There are many different embeddings of $SU(4)$ 
as a subgroup of $Spin(7)$, and there exist subgroups $SU(4)\subset 
Spin(7)$ which contain $H_2$ or $H_3$. Using these $SU(4)$ embeddings, 
we can construct 8-manifolds $Y_2$ and $Y_3$ desingularizing $\R^8/G$ 
in the same way as we made $Y_1$. At present, $Y_1,Y_2$ and $Y_3$ are 
the only ways the author knows to desingularize the singularity 
$\R^8/G$ within holonomy $Spin(7)$, and in particular it is unknown 
whether there exists any desingularization for which the na\"\i ve 
String Theory predictions hold.

The author believes that there is a generalization of the idea of
Weyl group of a quotient singularity, using which one can explain
the desingularizations $Y_1,Y_2$ and $Y_3$ in the following way,
using the method of \S 3. The Weyl group $W$ of the crepant 
resolution $X$ of $\C^4/H_1$ should be $\{1,\gamma\}\cong\Z_2$, where 
$\gamma$ multiplies by $-1$ in $H_2(X,\Z)$ and $H_6(X,\Z)$, and acts 
trivially on $H_4(X,\Z)$. The group homomorphism $\chi:K_1\rightarrow W$ 
must be $\chi(H_1)=1$, $\chi(\lambda H_1)=\gamma$, because this yields 
the correct action of $\lambda'$ on $H_*(X,\Z)$. We then see that 
the na\"\i ve String Theory predictions do not hold for $Y_1$ 
because they implicitly assume that $\chi\equiv 1$, but in fact 
$\chi$ is nontrivial.

\section{Orbifolds and string theory}

String Theory is a branch of high-energy theoretical physics in 
which particles are modelled not as points but as 1-dimensional 
objects -- `strings' -- propagating in some background space-time, 
which is usually curved and may have dimension 10, 11 or 26, 
depending on the theory. String theorists aim to construct a
{\it quantum theory} of the string's motion. The process of
quantization is extremely complicated, and fraught with
mathematical difficulties that are as yet still poorly understood.

String theorists believe that to each compact Calabi-Yau 3-fold
$X$ one can associate a {\it conformal field theory} (CFT), which 
is a Hilbert space with a collection of operators satisfying some 
relations, to be regarded as the quantum theory of strings moving 
in $X$. They can then use their understanding of conformal field 
theories to make conjectures about Calabi-Yau 3-folds, which have 
often turned out to be true, and mathematically very interesting.

String Theory can be used to study Calabi-Yau orbifolds, and their 
resolutions. The idea is this. Let $X$ be a compact Calabi-Yau
3-fold and ${\cal H}_X$ the associated CFT, and suppose that
$G$ is a finite group acting on $X$ preserving the Calabi-Yau
structure. Then $G$ also acts on ${\cal H}_X$, and Dixon et
al.~\cite{DHVW} showed that by a complicated process one
can construct a nonsingular quotient CFT ${\cal H}_{X/G}$ 
that corresponds to the singular Calabi-Yau orbifold~$X/G$. 

Then small, smooth deformations of ${\cal H}_{X/G}$ as a 
CFT correspond to desingularizations of $X/G$ as a Calabi-Yau
orbifold. Because of this, topological data such as the Hodge
numbers of these Calabi-Yau desingularizations of $X/G$ can
be extracted from ${\cal H}_{X/G}$, which in turn depends
only on $X$ and $G$. Therefore Dixon et al.~and others were
able to make predictions such as \eq{peuler} about the topology 
of Calabi-Yau desingularizations of~$X/G$.

Now it turns out that the quotient CFT ${\cal H}_{X/G}$ 
is not always uniquely determined. Instead, Vafa and Witten
\cite{Vafa2,VaWi} show that there can be several different
ways to reassemble the pieces of ${\cal H}_X$ to make 
${\cal H}_{X/G}$, a phenomenon which they call {\it discrete 
torsion}. One of these ways is preferred, corresponding to
crepant resolutions of $X/G$, but the other possibilities
correspond to certain partial desingularizations of $X/G$
in which a finite number of nodes remain unresolved, as
in~\S 5.

It would be interesting to be able to properly explain all
the different ways to desingularize a Calabi-Yau orbifold
in terms of String Theory, and hence perhaps to get a better 
understanding of these desingularizations. The outline of
this explanation already seems clear: to an orbifold $X/G$
we should be able to associate a finite number of different
quotient conformal field theories ${\cal H}_{X/G}$, and
every Calabi-Yau desingularization of $X/G$ should be 
associated to a small deformation of one of these theories.

However, there is a problem. The idea of discrete torsion 
gives a way to construct a small number of possibilities for 
${\cal H}_{X/G}$, but not enough different theories to 
explain the wide range of desingularizations possible in
examples. For instance, Vafa and Witten \cite[\S 2]{VaWi}
find only two possibilities for ${\cal H}_{X/G}$ in 
the case of $T^6/\Z_2^2$, but the analysis of \S 6
suggests that there should be hundreds of thousands of
orbifold theories. One explanation of this is given by Aspinwall 
and Morrison \cite[p.~128]{AsMo}. They observe that the CFT 
corresponding to the Calabi-Yau deformation of $T^6/\Z_2^2$ that 
they consider, is singular for the orbifold itself. Thus, one 
cannot construct a nonsingular orbifold CFT ${\cal H}_{X/G}$ 
corresponding to this family. Presumably this means that one must 
also consider {\it singular} CFT's~${\cal H}_{X/G}$.

I would like to finish by making two tentative suggestions
on how it may be possible to construct a larger set of
quotient conformal field theories ${\cal H}_{X/G}$
than is given by the framework of \cite{Vafa2,VaWi}. I am 
not competent to develop these ideas myself, but I hope 
that some physicist will do so and will tell me the answer. 
The first suggestion is very obvious: to regard discrete 
torsion as living locally on $X/G$ rather than globally 
in $H^2(G,U(1))$, so that one can make different choices 
of discrete torsion on different pieces of the singular 
set of~$X/G$.

The second suggestion is more serious: I think the ideas on 
Weyl groups in \S 3 should have an interpretation in
String Theory, and that this should lead directly to a way
to construct more quotient CFT's ${\cal H}_{X/G}$. In
particular, the Weyl group $W$ of a singularity should
appear in String Theory as an extra group of symmetries
acting on the quotient~CFT.

Suppose for instance that $G$ has a normal subgroup $H$ 
with $K=G/H$, and that the singularities of $X/H$ have Weyl 
group $W$. Then ${\cal H}_{X/H}$ should admit an action of 
$K\ltimes W$ preserving the CFT structure. Here $W$ acts nontrivially 
only on the twisted sectors in ${\cal H}_{X/H}$. The data
$\chi$ of \S 3 defines a subgroup $\tilde K$ of $K\ltimes W$
isomorphic to $K$, and quotienting ${\cal H}_{X/H}$ by $\tilde 
K$ should give an orbifold CFT corresponding to $X/G$. The difficult 
parts of this programme appear to be in understanding the action of
$W$ on ${\cal H}_{X/H}$, which is not obvious, and in making
sense of CFT sectors twisted by elements of $K\ltimes W$ rather than~$K$.

\end{document}